\documentclass[a4paper,english,10pt,oneside]{article}

\usepackage{amsmath}
\usepackage{amsthm}
\usepackage{amsfonts}
\usepackage{url}
\usepackage{color}
\usepackage{todonotes}
\usepackage{tabularx,supertabular,booktabs,array,longtable,multicol,multirow}

\usepackage{zibtitlepage}
\usepackage{geometry}
\usepackage[hidelinks]{hyperref}
\usepackage{caption}

\theoremstyle{definition}

\newtheorem{example}{Example}
\newtheorem{remark}{Remark}
\newcommand{\doilink}[1]{\href{http://dx.doi.org/#1}{{\textcolor{blue}{DOI:\allowbreak#1}}}}
\usepackage{xspace}
\newcommand{\solverfont}[1]{\scalebox{0.95}{#1}\xspace}
\newcommand{\scip}{\solverfont{SCIP}}
\newcommand{\scipc}{\solverfont{SCIP+C}}
\newcommand{\cplex}{\solverfont{CPLEX}}
\newcommand{\vipr}{\solverfont{VIPR}}
\newcommand{\lp}{\solverfont{LP}}
\newcommand{\mip}{\solverfont{MIP}}
\newcommand{\gmp}{\solverfont{GMP}}
\newcommand{\bmbf}{\solverfont{BMBF}}
\newcommand{\modal}{\solverfont{MODAL}}
\newcommand{\hollight}{\solverfont{HOL Light}}

\ZTPAuthor{Kevin K.~H.~Cheung, Ambros Gleixner, and Daniel E.~Steffy}
\ZTPTitle{\bf Verifying Integer Programming Results}

\ZTPInfo{This report has been published in the proceedings of IPCO
    2017. Please cite as:\\[1ex] K.K.H.~Cheung, A.~Gleixner, and D.E.~Steffy.
  Verifying Integer Programming Result.  In F.~Eisenbrand and J.~Koenemann,
  eds., {\em Integer Programming and Combinatorial Optimization: 19th
    International Conference, IPCO 2017}, pp.~148--160, 2017,
  \doilink{10.1007/978-3-319-59250-3\_13}}
\ZTPNumber{16-58}
\ZTPMonth{Revised version March}
\ZTPYear{2017}

\title{\Large\bf Verifying Integer Programming Results}

\author{\normalsize
  Kevin K.~H.~Cheung\thanks{School of Mathematics and Statistics, Carleton University, Ottawa, ON, Canada, \texttt{kevin.cheung@carleton.ca}.}
  \ and Ambros Gleixner\thanks{Department Optimization, Zuse Institute Berlin, Takustr.~7, 14195 Berlin, Germany, \texttt{gleixner@zib.de}, ORCID: \href{https://orcid.org/0000-0003-0391-5903}{0000-0003-0391-5903}.}
  \ and Daniel E.~Steffy\thanks{Department of Mathematics and Statistics, Oakland University, Rochester, MI, USA, \texttt{steffy@oakland.edu}.}
}

\date{\normalsize March 9, 2017}

\begin{document}

\zibtitlepage

\newgeometry{left=8.5pc,right=8.5pc,top=1.5in,textheight=640pt,footskip=17mm}

\maketitle

\begin{abstract}
Software for mixed-integer linear programming can return incorrect results 
for a number of reasons, one being the use of inexact floating-point arithmetic.
Even solvers that employ exact arithmetic may suffer from programming or 
algorithmic errors, motivating the desire for a way to produce independently 
verifiable certificates of claimed results.
Due to the complex nature of state-of-the-art \mip solution algorithms, the ideal form of such a
certificate is not entirely clear.
This paper proposes such a certificate format
designed with simplicity in mind, which is composed of a
list of statements that can be sequentially verified using a limited number of 
inference rules.
We present a supplementary verification tool for compressing and checking these
certificates independently of how they were created.  We report computational
results on a selection of \mip instances from the
literature.  To this end, we have extended the exact rational version of the \mip
solver \scip to produce such certificates.\\

\end{abstract}
 \section{Introduction}

The performance of algorithms for solving mixed-integer linear programs to
optimality has improved significantly over the last
decades~\cite{AchterbergWunderling2013,BixbyFenelonGuRothbergWunderling2000}.
As the complexity of the solvers increases, a question emerges:
\emph{How does one know if the computational results are correct?}

Although rarely, \mip solvers do occasionally return incorrect or dubious results~\cite{hybrid}.
Despite such errors, 
maintaining a skeptical attitude that borders on paranoia is arguably neither
healthy nor practical.  After all, machines do outperform humans on 
calculations by orders of magnitude and many tasks in life are 
now entrusted to automation.  Hence, the motivation for asking how to 
verify correctness of computational results is not necessarily because 
of an inherent distrust of solvers.  Rather, it is the desire to seek ways 
to identify and reduce errors and to improve confidence in the computed 
results.
Previous research on computing accurate solutions for \mip has utilized various 
techniques including interval arithmetic \cite{NS04}, 
exact rational arithmetic \cite{ExLP,hybrid,cert03}, and safely derived cuts 
\cite{safecuts}.
Nevertheless, as stated in \cite{hybrid},
``even with a very careful implementation and extensive testing,
a certain risk of an implementation error remains''.

One way to satisfy skeptics is formal code verification as is sometimes found in
software for medical applications and avionics.  For global optimization,
progress in this direction has been made very
recently~\cite{Narkawicz2013,Smith2015}.  For modern \mip solvers, which easily
consist of several 100,000~lines of code, this may be an ambitious goal.
An alternative is to build solvers that output extra information that
facilitates independent checking.  We 
shall use the word {\it certificate} to refer to 
such extra information for a given problem that has been solved.
Ideally, the certificate should allow for checking the results
using fewer resources than what are needed to solve the problem from scratch.  
Such a certificate could in principle be used in 
formal verification using a proof checker as done in
the Flyspeck Project~\cite{flyspeck,Obua,SoHa} 
for a formal proof of Kepler's Conjecture, or informal verification as done
by Applegate {\it et al.}~\cite{TSP} for the Traveling Salesman Problem 
and by Carr {\it et al.}~\cite{CGPP} in their unpublished work for
\mip in general.  Naturally, certificates should be as simple to verify as 
possible if they are to be convincing.

We highlight two specific applications where solution verification is 
desirable.
First, Achterberg~\cite{achterberg} presented \mip formulations for
circuit design verification problems, for which solvers have been shown to return incorrect 
results~\cite{hybrid}.
Second, Pulaj~\cite{Pulaj} has recently used \mip to 
settle open questions related to Frankl's conjecture.
Software developed in connection with this paper has been successfully 
used to generate and check certificates for \mip models coming from both of these 
applications.

For linear programming,
duality theory tells us that
an optimal primal solution and an optimal dual solution are sufficient
to facilitate effective verification of optimality.
In the case of checking infeasibility, a Farkas certificate will do.
Therefore, verifying \lp results, at least
in the case when exact rational arithmetic is used, is rather
straightforward.
However, the situation with \mip is drastically different.
From a theoretical perspective, 
even though some notions of duality for \mip have been formulated~\cite{GR07},
small (i.e. polynomial size) certificates for 
infeasibility or optimality may not even exist.
As a result, there are many forms that certificates could take: a 
branch-and-bound tree, a list of derived cutting planes, a superadditive dual 
function, or other possibilities for problems with special 
structures such as pure integer linear programming
and binary programming~\cite{BE,nullstellensatz,klabjan,La04}.
Which format would be preferred for certificate verification is not entirely 
clear, and in this paper we provide reasoning behind our choice.

From a software perspective, \mip result certification is also considerably
more complicated than \lp certification.
Even though most solvers adopt the branch-and-cut paradigm, 
they typically do not make the computed branch-and-bound tree or generated cuts 
readily available, and they may also utilize many other techniques 
including constraint propagation, conflict analysis, or reduced cost fixing.  
Thus, even if a solver did print out all information used to derive its 
solution, a verifier capable of interpreting such information would itself be 
highly complex, contradicting our desire for a simple verifier.
As a result, other than accepting the results of an exact solver such as 
\cite{hybrid}, the best that many people can do today to ``verify'' the results 
of a solver on a \mip instance is to solve the instance by several different 
solvers and check if the results match or minimally check that a returned
solution is indeed feasible and has the objective function value claimed,
as is done by the solution checker in \cite{miplib2010}.

The main contribution of this paper is the development of a certificate format 
for the verification of mixed-integer linear programs.
Compared to the previous work of Applegate {\it et al.}~\cite{TSP} for the 
Traveling Salesman Problem and the unpublished work of Carr {\it et 
al.}~\cite{CGPP} for general \mip, our certificate format has a significantly 
simpler structure.
It consists of a sequence of statements that can be verified one by one using 
simple inference rules, facilitating verification in a manner akin to natural 
deduction.
The approach is similar to that for verification of unsatisfiability proofs for 
SAT formulas. (See for example \cite{He13,We14}.)
This simple certificate structure makes it easier for researchers to 
develop their own independent certificate verification programs, or 
check the code of existing verifiers, even without any expert knowledge of \mip 
solution algorithms.

To demonstrate the utility of the proposed certificate format,
we have developed a reference checker in C++ and
added the capability to produce such
certificates to the exact version of the \mip solver \scip~\cite{hybrid,scip}.
We used these tools to verify results reported in \cite{hybrid}.
To the best of our knowledge, this work also represents the first software 
for general \mip certificate verification that has been made available 
to the mathematical optimization community.

{\it Organization of the paper.}
Even though the proposed format for the certificate is
straightforward, some of the details are nevertheless technical.
Therefore, in this paper, we discuss the certificate format
at a conceptual level.  The full technical specification is
found in the accompanying computer files.\footnote{See \url{https://github.com/ambros-gleixner/VIPR}}
We begin with the necessary
ingredients for the simple case of \lp in Section~\ref{lp}.
In Section~\ref{cg}, the ideas for dealing with \lp are extended 
to pure integer linear programming.
The full conceptual description of the format of the certificate
is then given in Section~\ref{milp}.
Computational experiments are reported in Section~\ref{compute},
and concluding remarks are given in Section~\ref{conclusion}.
Throughout this paper, we assume that problems are specified and solved
with exact rational arithmetic.
 \section{Certificates for linear programming}\label{lp}

A certificate of optimality for an \lp
is a dual feasible solution whose objective
function value matches the optimal value. 
However, there is no need to specify the dual
when one views the task of certification as an inference procedure,
see, e.g., \cite{Hooker}.
Suppose we are given the system of linear constraints
\[\tag{S}\label{equ:linsystem}
\begin{array}{c}
Ax  \geq  b,
A'x  \leq  b',
A''x  =   b'',
\end{array}
\]
where $x$ is a vector of variables,
\(A \in \mathbb{R}^{m\times n}\), 
\(A' \in \mathbb{R}^{m' \times n}\),
\(A'' \in \mathbb{R}^{m'' \times n}\),
\(b \in \mathbb{R}^{m}\), 
\(b' \in \mathbb{R}^{m'}\),
and \(b'' \in \mathbb{R}^{m''}\) for some nonnegative integers \(n\), \(m\),
\(m'\), and \(m''\).

We say that \(c^\mathsf{T} x \geq v\) is obtained by taking
a \emph{suitable linear combination} of the constraints in \eqref{equ:linsystem}
if \[c^\mathsf{T} = d^\mathsf{T} A + 
{d'}^\mathsf{T} A'
+ {d''}^\mathsf{T} A'',~
v = d^\mathsf{T} b + 
{d'}^\mathsf{T} b'
+ {d''}^\mathsf{T} b''\]
for some
\(d \in \mathbb{R}^m\),
\(d' \in \mathbb{R}^{m'}\), and
\(d'' \in \mathbb{R}^{m''}\) with
\(d \geq 0\) and \(d' \leq 0\).
If \(x\) satisfies \eqref{equ:linsystem}, then it necessarily satisfies
 \(c^\mathsf{T} x \geq v\).
We say that the inequality \(c^\mathsf{T} x \geq v\) 
is \emph{inferred} from \eqref{equ:linsystem}.  We will refer to this general inference procedure as \emph{linear inequality inference}.

\begin{remark}
Together, \(d,d',d''\) is simply a feasible solution to the linear 
programming dual of the linear program
\[\tag{LP}\label{equ:lp}
\min\{ c^\mathsf{T} x \mid Ax \geq b, A'x \leq b', A''x =  b'' \}.
\]
The inequality \(c^\mathsf{T} x \geq v\) is sometimes called a 
\emph{surrogate} of \eqref{equ:linsystem}.
(See \cite{Hooker}.)
\end{remark}

Suppose that an optimal solution to \eqref{equ:lp}
exists and the optimal value is~\(v\).
Linear programming duality theory guarantees that
\(c^\mathsf{T} x \geq v\)
can be inferred from \eqref{equ:linsystem}.
Therefore, linear inequality inference is sufficient to certify
optimality for linear programming.  
Conceptually, the certificate that we propose is a listing
of the constraints in \eqref{equ:linsystem} followed by the 
inequality \(c^\mathsf{T} x \geq v\) with the associated multipliers
used in the inference as illustrated in the following example.

\begin{example}\label{eg:lp}
The following shows an \lp problem and
its associated certificate.
\begin{center}
\begin{tabular}{ccc}
\(
\begin{array}{rll}
\min &  2x + y \\
\mbox{s.t.}  \\
  C1: & 5x-y \geq 2 \\
  C2: & 3x-2y \leq 1. \\
\end{array}
\)
&~~~~~~~~~~&\small
\begin{tabular}{|l|l|}\hline
\(
\begin{array}{rl}
{\bf Given } \\
  C1: & 5x-y \geq 2 \\
  C2: & 3x-2y \leq 1 \\
\end{array}
\) \\ \hline
\(
\begin{array}{rll}
{\bf Derived } && {\bf Reason}\\
  \texttt{obj}: & 2x+y \geq 1 & \{ 1\times C1 + (-1)\times C2\}
\end{array}
\) \\ \hline
\end{tabular}
\end{tabular}
\end{center}
Here, \(C1\) and \(C2\) are constraint labels.
Taking the suitable linear combination \(1 \times C1 + (-1)\times C2\) 
gives \(2x+y \geq 1\), thus establishing that \(1\) is a lower bound
for the optimal value.
\end{example}
\begin{remark}
This type of linear inference can also be used to derive 
$\leq$-inequalities or equality constraints.
Assuming that all problem data is rational, rational multipliers are 
sufficient to certify infeasibility or optimality.
\end{remark}

\section{Handling Chv\'{a}tal-Gomory cutting planes}\label{cg}
Gomory~\cite{Gomory} showed in theory that, for pure integer linear 
programming (IP), optimality or 
infeasibility can be established by a pure cutting-plane approach.
Such an approach can also work in practice~\cite{Balas,Zanette}.
In addition to linear inequality inference, a rounding operation
is needed.

Suppose that \(c^\mathsf{T} x \geq v\) can be inferred from \eqref{equ:linsystem}
by taking a suitable linear combination of the constraints.
If \(c_i \in \mathbb{Z}\) for \(i \in I\) for
some \(I \subseteq \{1,\ldots,n\}\) and 
\(c_i = 0\) for \(i \notin I\),
then any \(x \in \mathbb{R}^n\) satisfying \eqref{equ:linsystem} with
\(x_i \in \mathbb{Z}\) for \(i \in I\)
must also satisfy \(c^\mathsf{T} x \geq \lceil v \rceil\).
We say that \(c^\mathsf{T} x \geq \lceil v \rceil\) is obtained
from \(c^\mathsf{T} x \geq v\) by \emph{rounding}.
When \(I = \{1,\ldots,n\}\), the inequality 
\(c^\mathsf{T} x \geq \lceil v \rceil\)
is known as a \emph{Chv\'{a}tal-Gomory cut} (CG-cut in short).
It can then be added to the system and the process of obtaining another
CG-cut can be repeated.
Conceptually, a certificate for an IP instance solved using 
only CG-cuts can be given as a list of the original constraints
followed by the derived constraints.

\begin{example}\label{eg:ip}
The following shows an IP problem and
its associated certificate.
\begin{center}
\begin{tabular}{ccc}
\(
\begin{array}{rll}
\min & x + y \\
\mbox{s.t.}  \\
  C1: & 4x+y \geq 1 \\
  C2: & 4x-y \leq 2 \\
      & x,y \in \mathbb{Z}
\end{array}
\) 
& ~~~~~~~~~ &\small
\begin{tabular}{|l|l|}\hline
\(
\begin{array}{rl}
{\bf Given} \\
     & x,y \in \mathbb{Z} \\
 C1: & 4x+y \geq 1 \\
 C2: & 4x-y \leq 2 \\
\end{array}
\) \\ \hline
\(
\begin{array}{rll}
{\bf Derived } & & {\bf Reason}\\
 C3: & y \geq -\frac{1}{2} & \left\{ \frac{1}{2}\times C1 
                             + \left(-\frac{1}{2}\right)\times C2 \right \} \\
 C4: & y \geq 0 & \{ \mbox{round up } C3 \} \\
 C5: & x+y \geq \frac{1}{4} & \left\{ \frac{1}{4}\times C1 
                             + \frac{3}{4}\times C4 \right \} \\
 C6: & x+y \geq 1 & \{ \mbox{round up } C5 \} 
\end{array}
\) \\ \hline
\end{tabular}

\end{tabular} 
\end{center}

\end{example}

Note that the derived constraints in the certificate can be processed 
in a sequential manner.
In the next section, we see how to deal with
branching without sacrificing sequential processing.

\section{Branch-and-cut certificates}\label{milp}
In practice, most \mip instances are not solved by cutting planes 
alone.  Thus, certificates as described in the previous section are of
limited utility.
We now propose a type of certificate for optimality or infeasibility
established by a branch-and-cut procedure in which the generated
cuts at any node can be derived as split cuts and branching 
is performed on a disjunction of the form \(a^\mathsf{T} x \leq \delta
\vee a^\mathsf{T} x \geq \delta+1\) where \(\delta \in \mathbb{Z}\)
and \(a^\mathsf{T} x\) is integral for all feasible \(x\).

The use of split disjunctions allows us to consider branching and cutting 
under one umbrella.
Many of the well-known cuts generated by \mip solvers can be derived
as split cuts~\cite{corjols} and they are effective in closing
the integrality gap in practice~\cite{ricardo}.
Branching typically uses only simple split disjunctions (where 
the $a$ above is a unit vector), although some studies have considered 
the computational performance of branching on general disjunctions 
\cite{impgendisj,GamrathMelchioriBertholdGleixnerSalvagnin2015,gendisj,Ow01}.

Recall that each branching
splits the solution space into two subcases.
At the end of a branch-and-bound (or branch-and-cut) procedure,
each leaf of the branch-and-bound tree corresponds to one of the cases
and the leaves together cover all the cases that need to be considered.
Hence, 
if the 
branch-and-bound tree is valid, all one needs to look at are the \lp
results at the leaves.

Our proposal is to ``flatten''
the branch-and-bound tree into a list of statements that can be verified sequentially.
Thus, our approach departs from the approaches
in \cite{TSP} and  \cite{CGPP}, which require explicit
handling of the tree structure.  The price we pay 
is that we can no longer simply examine the leaves of the tree.
Instead, we process the nodes in a bottom-up fashion
and discharge assumptions as we move up towards the root.
We illustrate the ideas  with an example.

\begin{example}\label{cert_eg}
It is known that the following has no solution.
\[
\begin{array}{rr}
C1:&  2x_1 + 3x_2 \geq 1 \\
C2:&  3x_1 - 4x_2 \leq 2 \\
C3:&  -x_1 + 6x_2 \leq 3 \\
& x_1, x_2 \in \mathbb{Z}
\end{array}\]

Note that \((x_1,x_2) = (\frac{10}{17}, -\frac{1}{17})\) is an extreme point
of the region defined by \(C1\), \(C2\), and \(C3\).  Branching
on the integer variable  \(x_1\) leads to two cases: 

\begin{itemize}
\item {\bf Case 1.}  \(A1:~x_1 \leq 0\)

Note that \((x_1,x_2) = (0, \frac{1}{3})\) satisfies \(C1, C2, C3, A1\).
We branch on \(x_2\):

\begin{quote}
{\bf Case 1a.} \(A3: ~x_2 \leq 0 \)

Taking \( C1 + (-2)\times A1 + (-3)\times A3 \) gives
the absurdity \(C4:~0 \geq 1.\)

{\bf Case 1b.} \(A4: ~x_2 \geq 1 \)

Taking 
\( \left(-\frac{1}{3}\right) \times C3 + 
\left(-\frac{1}{3}\right)\times A1 + 2\times A4 \) 
gives the absurdity \(C5:~0 \geq 1.\) 
\end{quote}

\item {\bf Case 2.}  \(A2:~x_1 \geq 1\)

Taking 
\( \left(-\frac{1}{4}\right)\times C2 
+ \left(\frac{3}{4}\right)\times A2\)
gives \(C6:~ x_2 \geq \frac{1}{4}.\)
Rounding gives
\(C7:~ x_2 \geq 1.\)

Taking 
\( \left(-\frac{1}{3}\right)\times C2 
+ (-1)\times C3 + \frac{14}{3}\times C7\)
gives the absurdity \(C8:~ 0 \geq 1.\)
\end{itemize}
As all cases lead to \(0 \geq 1\), we conclude that there
is no solution.  To issue
a certificate as a list of derived constraints, we need a way to 
specify the different cases.  To this end, we allow the introduction
of constraints as assumptions.

Figure~\ref{fig:complete} shows a conceptual certificate for the instance.
Notice how the constraints \(A1\), \(A2\), \(A3\), and \(A4\) are introduced
to the certificate as assumptions.
Since we want to end with \(0 \geq 1\) without additional assumptions attached,
we get there by gradually undoing the case-splitting operations. 
We call the undoing operation \emph{unsplitting}.  
For example, \(C4\) and \(C5\) are both the absurdity
\(0 \geq 1\) with a common assumption \(A1\).  Since \(A3\vee A4\) is
true for all feasible \(x\), 
we can infer the absurdity \(C9:~0\geq 1\) assuming only
\(A1\) in addition to the original constraints. We say that
\(C9\) is obtained by \emph{unsplitting} \(C4,C5\) on \(A3,A4\).
Similarly, both \(C8\) and \(C9\) are the absurdity \(0 \geq 1\) and 
\(A2 \vee A1 \) is true for all feasible \(x\), we can therefore
unsplit on \(C8,C9\) on \(A2,A1\) to obtain \(C10:~0 \geq 1\) 
without any assumption in addition to the original constraints.  

\begin{figure}[htbp]
\begin{center}\small
\begin{tabular}{|l|l|}\hline
\(
\begin{array}{rl}
{\bf Given} \\
    & x,y \in \mathbb{Z} \\
 C1:&  2x_1 + 3x_2 \geq 1 \\
 C2:&  3x_1 - 4x_2 \leq 2 \\
 C3:&  -x_1 + 6x_2 \leq 3 \\
\end{array}
\) \\ \hline
\(
\begin{array}{rlll}
{\bf Derived} & & {\bf Reason} & {\bf Assumptions}\\
 A1: & x_1 \leq 0 & \{ \mbox{assume} \} \\
 A2: & x_1 \geq 1 & \{ \mbox{assume}\} \\
 A3: & x_2 \leq 0 & \{ \mbox{assume}\} \\
 C4: & 0 \geq 1 & \left\{ C1+ (-2)\times A1 + (-3)\times A3 \right \} 
                    &  A1, A3 \\
 A4: & x_2 \geq 1 & \{ \mbox{assume}\} \\
 C5: & 0 \geq 1 & \left\{ \left(-\frac{1}{3}\right) \times C3 + 
                  \left(-\frac{1}{3}\right)\times A1 + 2\times A4 \right \} 
                    &  A1, A4 \\
 C6: & x_2 \geq \frac{1}{4} & \left \{ \left(-\frac{1}{4}\right)\times C2 
                           + \left(\frac{3}{4}\right)\times A2 \right \} 
                    &  A2  \\
 C7: & x_2 \geq 1 & \left \{ \mbox{round up } C6 \right \} & A2 \\
 C8: & 0 \geq 1 & \left \{ \left(-\frac{1}{3}\right)\times C2 
                         + (-1)\times C3 + \frac{14}{3}\times C7
                   \right \} &  A2  \\
 C9: & 0 \geq 1 & \left \{ \mbox{unsplit } C4, C5 \mbox{ on }
A3, A4 \right \} &  A1 \\
 C10: & 0 \geq 1 & \left \{ \mbox{unsplit } C8, C9 \mbox{ on }
A2, A1 \right \} \\
\end{array}
\) \\ \hline
\end{tabular}
\end{center}
\caption{Certificate for Example~\ref{cert_eg}}
\label{fig:complete}
\end{figure}

\end{example}

In practice, the list of assumptions associated with each derived constraint 
needs not be specified explicitly in the certificate, but can be deduced on the fly
by a checker.
For example, when processing \(C4\), we see that it uses \(A1\) and
\(A3\), both of which are assumptions.
Hence, we associate \(C4\) with the list of assumptions \(A1,A3\).
As any linear inequality can be introduced as an assumption, 
branching can be performed on general disjunctions.

\begin{remark}
  Our proposed certificate can also be used to represent split cuts.
  Split cuts are inequalities that are valid for the defining inequalities
  taken together with each one of the inequalities in a split disjunction,
  $ a^\mathsf{T} x \leq \delta \vee a^\mathsf{T} x \geq \delta+1$, where $\delta$ is an integer and $a$
  is an integer vector that is nonzero only in components corresponding to integer variables.
  To derive a proof of a split cut's validity, the inequalities in the split
  disjunction can each be introduced as assumptions, the cut can be derived
  for each side of the split disjunction using linear inequality inference,
  and then unsplitting can be applied to discharge the assumptions.
\end{remark}

 \section{Computational experiments}\label{compute}

In this section, we describe software developed to produce and check 
certificates for \mip results using the certificate format 
developed in this paper.
It is freely available for download, 
along with a precise technical specification of
the file format.\footnote{See \url{https://github.com/ambros-gleixner/VIPR}}
One of its features is that after each derived 
constraint an integer is printed to specify the largest index of any derived constraint that references it.  This allows constraints to be freed from memory 
when they will no longer be needed. The following C++ programs are provided:

\begin{itemize}
\item {\tt viprchk} verifies \mip results 
provided in our specified file format.  All computations are performed in exact 
rational arithmetic using the \gmp library~\cite{gmp}.
\item {\tt viprttn} performs simple modifications to ``tighten'' certificates.
  It removes unnecessary derived constraints to reduce the file size.  In order
  to decrease peak memory usage during checking, it reorders the remaining ones
  using a depth-first topological sort and for each derived constraint that
  remains, it computes the largest index over constraints that references it.
\item {\tt vipr2html} converts a certificate file to a 
 ``human-readable'' HTML file.
\end{itemize}

We again emphasize that the format was designed with simplicity in mind; 
the certificate verification program we have provided is merely a reference 
and others should be able to write their own verifiers without much difficulty.

In addition, we created a modified version of the exact rational \mip solver 
described in \cite{hybrid} and used it to compute certificates for several 
\mip instances from the literature.
The exact rational \mip solver is based on \scip~\cite{scip} and uses a 
hybrid of floating-point and exact rational arithmetic to efficiently compute 
exact solutions using a pure branch-and-bound algorithm.
In our experiments, the rational \mip solver uses \cplex~12.6.0.0
\cite{cplex} as its 
underlying floating-point \lp solver and a modified version of 
QSopt\_ex~2.5.10
\cite{QSopt}  as its underlying exact \lp solver.
The exact \mip solver supports several methods for computing valid dual bounds 
and our certificate printing functionality is currently supported by the 
\emph{Project-and-shift} method (for dual solutions only) and the \emph{Exact 
LP} method (for both dual solutions and Farkas proofs), for details 
on these methods see \cite{hybrid}.
This developmental version is currently available from the authors by 
request.
We note that the certificate is printed concurrently with the solution process
which leads to certificates that have potential for reduction and
simplification by {\tt viprttn}, or other routines.
For example, as each node is processed its derived dual bound is printed to
the certificate even though it may become redundant if branching is performed
and new dual bounds are computed at the child nodes; also, discovery of a new
primal solution might allow pruning of a large subtree, rendering many bound
derivations redundant.

The program {\tt viprttn} processes the list of derived constraints in two
passes.  In the first pass, it builds the dependency graph 
with nodes representing the derived constraints and arcs
\(uv\) such that the derived constraint represented by
\(u\) is referenced by the reason for deriving the constraint represented
by \(v\).  In the second pass, it performs a topological
sort using depth-first search on the component that contains the final 
constraint and writes out the reordered list of derived constraints
with updated constraint indices.  

\newcommand{\numinst}{\ensuremath{N}}
\newcommand{\numsolved}{\ensuremath{N_\text{sol}}}
\newcommand{\sciptime}{\ensuremath{t_{\text{MIP}}}}
\newcommand{\ttntime}{\ensuremath{t_\text{ttn}}}
\newcommand{\chktime}{\ensuremath{t_\text{chk}}}
\newcommand{\origsize}{\ensuremath{\text{size}_\text{raw}}}
\newcommand{\ttnsize}{\ensuremath{\text{size}_\text{ttn}}}
\newcommand{\gzsize}{\ensuremath{\text{size}_\text{gz}}}
In the following, we report some computational results on the time and 
memory required to produce and verify certificates.
We considered the \emph{easy} and \emph{numerically difficult} (referred to 
here as `\emph{hard}') test sets from \cite{hybrid}; these test sets 
consist of instances from well known libraries including
\cite{miplib2003,miplib3,miplib2010,coral,mittelmann}.
Experiments were conducted on a cluster of Intel(R) Xeon(R) CPU E5-2660 v3 at 
2.60\,GHz; jobs were run exclusively to ensure accurate time measurement.
Table~\ref{tab:results} reports a number of aggregate statistics on these 
experiments. 
The columns under the heading \scip report results from tests using the exact 
version of \scip, using its default dual bounding strategy.
The columns under \scipc report on tests involving the version of exact \scip 
that generates certificates as it solves instances;
it uses only the dual bounding methods \emph{Project-and-shift} and \emph{Exact LP}
that support certificate printing, contributing to its slower 
speed.
Columns under the heading \vipr report time and memory usage for certificate 
checking.

For each of the easy and hard test sets, we report information aggregated into 
four categories: `all' reports statistics over all instances; `solved' reports 
over instances solved by both \scip and \scipc within a 1~hour time limit and a 
10\,{\small GB} limit on certificate file size; `memout' reports on 
instances 
where \scipc stopped because the certificate file size limit was reached;
and `timeout' reports on the remaining instances, where one of the solvers hit the time limit.
All averages are reported as shifted geometric means with a shift of 10~sec.\ 
for time and 1\,MB for memory. 
The column \numinst{} represents the number of instances in each category; 
\numsolved{} represents the number in each category that were solved to 
optimality (or infeasibility) by a given solver; 
\sciptime{} represents the time (sec.) used to solve the instance and, when 
applicable, output a certificate; 
\ttntime{} is the time (sec.) required by the \texttt{viprttn} routine to 
tighten the certificate file; 
\chktime{} is the time (sec.) required to for \texttt{viprchk} to check the 
certificate file -- on instances in the memout and timeout rows this represents 
the time to verify the primal and dual bounds present in the intermediate 
certificate printed before the solver was halted.
The final three columns list the size of the certificate (in MB), before 
tightening, after tightening and then after being compressed to a gzipped file.
Timings and memory usage for individual instances are available 
in a document hosted together with the accompanying software.

\begin{table}
  \caption{Aggregated computational results over 107~instances from~\cite{hybrid}.}
  \label{tab:results}
  \setlength{\tabcolsep}{1pt}

  \small
  \begin{tabular*}{\textwidth}{@{\extracolsep{\fill}}lrrrrrrrrrr}
    \toprule
    & & \multicolumn{2}{c}{\scip} & \multicolumn{2}{c}{\scipc} & \multicolumn{5}{c}{\vipr}\\
    \cmidrule(){3-4}\cmidrule(){5-6}\cmidrule(){7-11}
    Test set               & \numinst & \numsolved & \sciptime & \numsolved & \sciptime & \ttntime & \chktime & \origsize & \ttnsize & \gzsize\\
    \midrule
    easy-all               &       57 &         54 &      63.3 &         39 &     190.9 &      8.9 &     27.2 &       227 &       77 &      24\\
    \phantom{easy}-solved  &       39 &         39 &      23.2 &         39 &      48.0 &      3.6 &     11.5 &        77 &       34 &      10\\
    \phantom{easy}-memout  &        5 &          4 &     600.6 &          0 &    1760.4 &     47.8 &    138.3 &     10286 &      513 &     157\\
    \phantom{easy}-timeout &       13 &         11 &     338.3 &          0 &    3600.0 &     23.3 &    102.7 &      1309 &      434 &     129\\
    \midrule
    hard-all               &       50 &         23 &     725.2 &         14 &     975.6 &      7.9 &     12.1 &       373 &       38 &      11\\
    \phantom{hard}-solved  &       13 &         13 &      22.9 &         13 &      40.7 &      2.2 &      5.3 &        49 &       15 &       5\\
    \phantom{hard}-memout  &       12 &          2 &    2476.4 &          0 &    1713.1 &     32.7 &     59.8 &     10266 &      235 &      67\\
    \phantom{hard}-timeout &       25 &          8 &    2052.1 &          1 &    3518.1 &      4.3 &      5.3 &       216 &       25 &       7\\
    \bottomrule
  \end{tabular*}
\end{table}

From this table, we can make a number of observations.  First, there 
is a noticeable, but not prohibitive, cost to generate the certificates.
The differences in \sciptime{} between \scip and \scipc are due to both the 
difference in dual bounding strategies, and the overhead for writing the 
certificate files.
In some additional experiments, we observed that on the 39 instances in the 
easy-solved category, the file I/O amounted to roughly 7\% of the solution 
time, based on this we believe that future modifications to the code will allow 
us to solve and print certificates in times much closer to those in the \scip 
column.
Perhaps most importantly, we observe that the time to check the certificates is 
significantly less than the time to solve the instances.

Moreover, the
certificate tightening program \texttt{viprttn} is able to make significant 
reductions in the certificate size, and the resulting certificate sizes are 
often surprisingly manageable.
Most striking is the tightening in the memout categories, which significantly
exceed the approximately 50\% reduction that could be expected by removing the
redundant linear inferences derived for internal nodes of the branch-and-bound
tree.
The most extreme tightening was achieved for the instance \texttt{markshare1\_1}
in `easy-memout', from 10\,{\small GB} to 8\,{\small kB}.
This is explained by the fact that the root dual bound is already zero and the
tree search is only performed for finding the optimal solution.  Hence, the
certificate is highly redundant and the derived constraints for all but the root
node can be removed.

The average reductions in the other categories are smaller, but also strictly
above 50\%.  This shows that \texttt{viprttn} performs more than just a removal of
internal nodes.
These results also show two aspects in which \scip's certificate printing
can be improved: by avoiding printing dual bound derivations for internal nodes
using a buffering scheme, and by not generating dual bound derivations for nodes
that do not improve upon the bound of the parent node.
 \section{Conclusion}\label{conclusion}
This paper presented a certificate format for verifying integer 
programming results.  We have demonstrated the practical feasibility of 
generating and checking such certificates on well-known \mip instances.  
We see this as the first step of many in verifying the results of integer 
programming solvers.
We now discuss some future directions made possible by this work.

Even in the context of floating-point arithmetic, our 
certificate format could serve a number of purposes.
Using methods described by \cite{safecuts,NS04}, directed rounding and interval 
arithmetic may allow us to compute and represent valid certificates exclusively 
using floating-point data, allowing for faster computation and smaller 
certificate size.
Additionally, generating approximate certificates with inexact data could be 
used for debugging solvers, or measuring the maximum or average 
numerical violation over all derivations.
In a more rigorous direction, one could also convert our certificates to a form 
that could be formally verified by a proof assistant such as \hollight~\cite{hollight}.

\paragraph*{Acknowledgements.}
We thank Kati Wolter for the exact version of \scip~\cite{hybrid} and
Daniel Espinoza for {QSopt\_ex}~\cite{QSopt}, which provided the basis for
our experiments, and Gregor Hendel for his Ipet
package~\cite{Hendel2014}, which was a big help in analyzing the experimental
results.
This work has been supported by the Research Campus \modal \emph{Mathematical
  Optimization and Data Analysis Laboratories} funded by the Federal Ministry of
Education and Research (\bmbf Grant~05M14ZAM).  All responsibility for the content
of this publication is assumed by the authors.

\end{document}